\documentclass[10pt]{article}
\usepackage{amssymb,amsfonts}
\usepackage{pb-diagram}
\usepackage{amsmath,amsthm,amsfonts,amssymb}
\usepackage[usenames]{color}
\usepackage{hyperref}
\usepackage{soul}
\usepackage{graphicx}

\newcommand{\be}{\begin{enumerate}}
\newcommand{\ee}{\end{enumerate}}
\newcommand{\beq}{\begin{equation}}
\newcommand{\eeq}{\end{equation}}
\usepackage{amssymb}
\usepackage{amsmath}
\usepackage{amssymb}

\newcommand{\comment}[1]{} % for multi-line comments
\begin{document}

\title{On Tarski's Decidability Problem:  Response to Sela's Report}
\author{Olga Kharlampovich, Alexei Myasnikov}

\maketitle

\begin{abstract}
This note provides a brief guide to the current state of the literature on  Tarski's problems with emphasis on features that distinguish the approach based on combinatorial  and algorithmic group theory from the topological approach to Tarski's problem. We use this note to provide corrections to some typos and to address some misconceptions from the recent report by Z. Sela about the relations between the concepts and results in the approaches to the Tarski problems. We were forced to read Sela's papers to be able to address some of his comments, and found  errors in his papers 6, 3 and 4 on Diophantine Geometry published in GAFA and  Israel J. Math. which we mention in Section 4. 
His proceedings of the ICM 2002  paper also contains wrong Theorem 6 (to make it correct one has to change the definition of non-elementary hyperbolic $\omega$-residually free towers to make them equivalent to our coordinate groups of regular NTQ systems.)
\end{abstract}
\section{Introduction} 

This is a response to Sela's  ``report'' \cite{Selar} on our paper ``Elementary theory of a free group'' \cite{KMel}.

Sela claims that the paper contains 57 mistakes that invalidate our results. These claims in some cases address inaccuracies in formulations (1, 5-8, 11 etc.) that were clarified  later in the text or in the subsequent literature (\cite{BKM}, \cite{KMqe}) or are misprints (2,3,12 etc.), some corrected later \cite{KMqe}; in many cases they stem from misunderstanding or deliberate misreading the conditions (9-10, 13-15 etc.). 

We discuss all  these points in the following section.   We are sure that Sela may immediately come up with another list of  ``mistakes'' of the same kind. This has already been going on for about 10 years and has very little to do with mathematics.

In Section 3 we discuss the state of the theory of Algebraic Geometry in free groups from our point of view. 

Finally we mention some serious inaccuracies and ``similarities'' in Sela' s papers that we are aware of. Possibly they can be fixed. We just want to say that ``those who live in glass houses should not throw stones''.

\section{About Sela's comments in \cite{Selar}}

 We emphasize the minor nature of the corrections below; any unbiased reader can verify that these in no way detract from the correctness of our main theorem. We keep the same numeration of comments as in \cite{Selar}.

(1) p. 470, Lemma 8.  
The lemma is formulated for $F$-groups and $F$-automorphisms, namely automorphisms fixing elements of the subgroup $F$ that is the free group of coefficients, therefore the only inner $F$-automorphism is the identical automorphism.
The slight problem in the lemma was corrected in our paper \cite{BKM} (see part of Theorem 5.3  in \cite{BKM}).  The error is completely inessential, we proved correctly that subgroups in question commute, which is the only result from this lemma that was used in our proof. The minor inaccuracy concerns  the intersection of the subgroups, it was never used or mentioned later in the paper.

(2) p. 482, Lemma 18. There is, indeed, a misprint in the formulation of Lemma 18.   In line -1 
one has to add "non--transcendental".  It is clear from p. 483, lines 8-10 that in the case of transcendental $\kappa$ we  consider the quotient of $K$, and $K_1$ is changed too.  In the rest of the paper the statement of this lemma is only used in the proof of Lemma 19 for non-transcendental solutions. A careful reader would be able to correct the formulation himself/herself looking at what is written in the proof.

(3) p. 483, Lemma 19. The misprint in the Definition 20 (of reducing quotients) was noticed a long time ago by N. Touikan and already corrected, see Definition 1 in \cite{KMqe}, where item (3) from Definition 20 was removed. This is not essential for the understanding of the proof.

(4) p. 485 lines -5 to -3. `` since the group $K$ is $\ldots$ one may assume that $\ldots$ is a monomorphism.''  There is no problem here. We do not say that the homomorphism $\kappa$ is a monomorphism.
One may assume this in the proof because in the case when it is not we either have a reducing quotient of $K$ and the theorem is proved, or we just consider $\kappa (K)$ instead of $K$. Sela's comment here is the result of misinterpretation. Our argument here is quite clear.

(5) Since comment (4) is dismissed,  we will talk only about the effectiveness of the bound. We agree that the proof of the effectiveness of the bound is rather sketchy because it requires good understanding of paper \cite{Imp}.  On the other hand, we did not provide the detailed proof because  the fact that the bound on the number of  the Max equivalence classes can be found follows from the existence of this bound and the decidability of the AE-theory of a free group (even though the decidability of the AE-theory is proved later in the paper, the algorithm for  computing this bound is not used in the proof of the decidability of the AE-theory).  We were asked by different people to explain this in detail, and the explanation is given in \cite{KMqe}. 

(6)-(8) The accusation here is that we use the term `` induced NTQ system'' and Sela has the term `` induced resolution''.  Our term came from the Bass-Serre theory,  the subgroups of fundamental groups of graphs of groups have induced decomposition, the notion has been in use since 1970. Notice that using his logic one can accuse him in using the term ``Zarisski closure'' without any references to the paper \cite{BMR} where it was shown that the Zarisski topology can be defined on $F^n$.

  Sela overlooked that it was said before (page 504, line -11 and page 505, line 9) that we would only consider fundamental sequences  satisfying 1'st and 2'nd restrictions.  
The phrase
``We then add conjugating elements so that all rigid subgroups that are mapped to the same free factor on the next level are conjugate into one subgroup'' is indeed missing in our description  on page 506, but  we say on page 507, line 9 that we end up with an NTQ system, which means  that such conjugating elements must be added (unless one is trying to present inaccuracies as ``mistakes''). 

(9)-(10) There is no problem here.
We say that we consider levels from the top to the bottom, obviously we are talking about the decrease of some characteristic on some level from the top, and on the higher levels the decompositions are the same as on the previous step. The comment is just a misinterpretation.

(11)    Once again, Sela overlooked that it was said before (page 504, line -11 and page 505, line 9) that we would only consider fundamental sequences  satisfying 1'st and 2'nd restrictions.  

(12)  This is indeed a regrettable misprint, we  consider all the free factors together and not continue with each factor separately. But everything that is said immediately after line 22  is correct, we can construct finitely many NTQ systems as in the paper etc. So the wrong statement is never used. 
Moreover everywhere in the proof we consider these free factors together.

(13)-(15)  Lemma 23,  Everything is correct here. 
There is the following  idea behind this. Let a solution $X$ of the equation $U_0(X)=1$ factors through a fundamental sequence (satisfying 1'st and 2'nd restrictions). In general we cannot lift a solution of $U(X,Y)=1$  into the group $G=<X|R(U_0)>$.   We can only obtain a formula solution in a corrective extension of the NTQ group corresponding to this fundamental sequence (by our parametrization theorem). Therefore we have to add extra  variables to $X$. Each specialization $\bar X$ of $X$ that factors through the fundamental sequence above can be extended to different specializations of these extra variables.  But all these different specializations form a finite number of fundamental sequences modulo $G$ (or, equivalently, we can consider fundamental sequences for $<X_2,\ldots, X_k>$ modulo rigid subgroups in the JSJ decomposition of $G$ and edge groups connecting rigid and abelian vertex groups). These fundamental sequences (denote them $b_1,\ldots,b_s$) end with groups with no sufficient splitting. Equation $U_1(X_1,\ldots, X_k)=1$ on page 520, line -3-5 is obtained from the following condition: ``we collect ALL values of $Y$ given by the formulas ...", such that for any specialization of  fundamental sequences $b_1,\ldots,b_s$ and any formal solution $Y=f(X_1,\ldots,X_k)$ depending on this specialization (that also belong to a finite number of fundamental sequences) we have equation $V,(X_1,\ldots,X_k, f(X_1,\ldots,X_k))=1$  By the parametrization theorem  (used twice here)  we have equation $U_1(X_1,\ldots ,X_k)=1$ and an equation on terminal groups of fundamental sequences $b_1,\ldots ,b_s$ (the terminal groups of their reducing quotients are included here).

The proof of the lemma shows that in the case when $H_t$ is not a proper quotient of $G$, we can replace $H_t$ be a finite number of quotients such that rigid subgroups of the terminal groups of $b_1,\ldots ,b_s$ are elliptic in the JSJ decompositions of these quotients (if not, this means that the homomorphisms factor through a reducing quotient of the terminal group of the corresponding $b_i$). Then we can apply the automorphisms corresponding to the decomposition $D_t$, and this is equivalent to taking minimal homomorphisms corresponding to  terminal groups  of fundamental sequences $b_1,\ldots ,b_s$. We are not loosing values of the initial variables $X$ that are candidates for making the sentence $\Phi$ true.

This is more complicated than what Sela advises, namely to begin this step by only considering  instead of $V_{fund} (U_1)$ minimal homomorphisms with respect to the fundamental sequences  modulo (free non-cyclic factors of) the image of $G$ in the corrective extension. It is indeed  sufficient to extend $\bar X$ to a minimal specialization of these extra variables, to get a formula solution, substitute it into $V(X,Y)=1$ and get an equation on this specialization of the variables. 
Nevertheless, what we are doing is not an error, and the counterexample is irrelevant.  Our restriction $U_1(X_1,\ldots, X_k)=1$ in general may be stronger than the one used by Sela.
Moreover, our construction is working in the general case with parameters (tree $T_{X_k}$ in Section 12.2) and we do not repeat it in the general case.

We do not understand the example. In the example, $L$ is itself an NTQ group, we would not extend it as Sela does. The maximal standard quotient is a free group
generated by $b_1,b_2$. We do not exactly know how to construct Sela's completion in this case (we could not find the explanation how it is constructed for that case in \cite{Sela2}), but $L\rightarrow \pi _1(S_2)$ is not even a fundamental sequence in our terminology.

(17)-(25)  These are all derivatives of the missing sentence in the construction of the induced NTQ system in comments (6)-(8).   In comment (20), about Case 3, there is a misinterpretation of what is done in the case when $N_0^1$ is mapped to a proper quotient. We consider fundamental sequences for this quotient as we did on step 1 for the proper quotient of the NTQ group. Namely, fundamental sequences for this quotient modulo the images of free factors of the subgroup generated by the images of $X_2,\ldots ,X_n$, and apply to them step 1. So these quotients are treated the same way as we treated proper quotients of NTQ groups, but not the same was as we treated proper quotients of $G$. 

Notice, that induced NTQ system is just a technical tool to prove that one can construct formula solutions over block-NTQ groups. 

(22) There is no any error. One, indeed, needs to look at the structure of the fundamental sequence induced by $G_{corr}$, and the "complexity" drops or the procedure stops along this branch. That's exactly what is written or implied here.  General comment: what we are doing in the construction of the $\forall\exists$-tree is the same as what we are doing in the construction of $T_{X_k}$, a similar tree but with parameters.  Sela is complaining that our process is not the same as his. 
This does not mean that our process does not converge. It looks like in $\forall\exists$-case Sela's process is simpler, but since our process is the same in  the  $\forall\exists$ case and in the case with parameters we don't have to write it twice.

(26)-(27)   It is clearly written in our paper: `` Consider a finite family of terminal groups of fundamental sequences for P modulo factors in the free decomposition of  $$F_{R(U)}=F_{R(\bar U_1)}\ast\ldots \ast F_{R(\bar U_k)} \ast F(Z),$$ where $F(Z)$ is a free group with basis $Z$, and $F_{R(\bar U_1)},\ldots ,F_{R(\bar U_k)}$ are freely indecomposable modulo...''.  Sela misinterprets saying  that the decomposition is modulo $F_{R(\bar U_1)},\ldots , F_{R(\bar U_k)}.$ In his example the only factor is $F(Z)$ and the constructed decomposition  is not modulo $F(Z)$.  

The free factors, indeed, can be conjugated, but the first and second restrictions on fundamental sequences for $F_{R(U)}$ are always assumed, see our comment on (6)-(8),
therefore they are conjugated into different free factors. 

(28) Our explanation here is unclear but  the statement that ``it does not increase the dimension'' is correct.

(29) Sela's critics here is that  constructing our tight enveloping fundamental sequence we work from  top to bottom, and he goes from bottom to top constructing his core resolutions.
First of all, several constructions don't work if to go from bottom to top. Secondly,
we repeat our construction several times (=several iterations) unlit the dimension stabilizes (see p. 529, lines -15--14),  this gives a fundamental sequence that has the desired properties. That might be similar to doing some constructions from bottom to top, because extra elements added to  level $j$  on iteration number $i$ are taken into account considering level $j+1$ on iteration $i+1$.

(30)-(31) First and second restrictions on fundamental sequences are always assumed, see our comment on (6)-(8). There is a particular reference to them again in lines -4-3, page 529.

(32)-(33) Derivatives of the misunderstanding of our induced fundamental sequence.
Our tight enveloping fundamental sequences, although constructed differently, satisfy the first two properties of Sela's core resolution. It is hard to say if they satisfy property (iii) because in property (iii), Definition 4.1
\cite{Sela5} he considers some completion $Comp(Res)(t,v,a)$ that has not appear before in the definition, and it is not the completion of $Res(t,v,a)$ because that one is denoted $Comp(Res)(u, t,v,a)$, this property cannot be understood.

(34)  Our tight enveloping fundamental sequences are firm in Sela's definition. They are similar to  his "penetrated" core resolutions.  At least they are as good for the purpose of proving that the process stops.  

(35)  The only possible meaning of the minimality  of solutions modulo some group  (in this case modulo $F_{R(U_i)}$) when we work with parameters, can be  that the corresponding variables belong to the terminal group of a fundamental sequence modulo  $F_{R(U_i)}$ (this group does not have a sufficient splitting but may have a splitting) and not factoring through a reducing quotient of this terminal group. So when we say that some solutions are minimal modulo $F_{R(U_i)}$ , that's what we mean.   

Indeed, we forgot to write on page 529 , line 8 about the possible reducing quotients, the sentence there should be:
 The group
$F_{R(M_i)}$ does not split either as a free product or as a
centralizer extension  modulo $F_{R(\bar U_i)}$, and if it has a
splitting, {\bf we consider fundamental sequences modulo $F_{R(\bar U_i)}$, their terminal groups and } only solutions of $F_{R(M_i)}$ minimal with
respect to the groups of canonical automorphisms corresponding to the fundamental sequence and its
 terminal group if it has a splitting {\bf and not factoring through a reducing quotient of this terminal group}. 
 
 As  we understand Sela calls these terminal groups Auxiliary limit groups (but we do not understand for what purpose he needs ``auxiliary resolutions'' unless they are one-step resolutions of the auxiliary limit group).

(37)-(46) Derived from the misunderstanding of our  tight enveloping fundamental sequences.

(43) We change the tight enveloping fundamental sequences  a bit  when compare $TEnv(W_{n-1})$ and $TEnv(W_{n})$ on page 531. At this point they resemble more Sela's ``core resolutions'' than ``penetrated core resolutions''.
Notice that he  defines both notions  many times and every time in different way, so it is hard to say what it is exactly.

(47) All what is said on page 541 line 18 is that there is a generic family of solutions
for the enveloping system. And there is nothing wrong with this.  We do not know why Sela calls our implicit function theorem or parametrization theorem a general form of Merzlyakov's theorem. He definitely prefers to refer to Merzlyakov who did not prove the theorem rather than for us who proved it.  Merzlyakov's theorem is about positve sentences.   We know how to generalize this theorem to the case of ``reduced modular groups''. We  do not introduce the notion of ``framed'' fundamental sequences and framed generic families because this is a technicality.  But when we lift solutions in block-NTQ groups (see page 522, line 16) we implicitly use similar techniques. (Notice also, that  the number of possible ``frames'' is finite as well as the number of ``framed'' fundamental sequences on each step).  

It worth mentioning that Sela makes a mistake in \cite{Sela4}, page 95, lines -8- (-4).   He states that we can associate the set of closures and formal solutions with the resolution CRes  and refers to Theorem 1.18 in \cite{Sela2}. This theorem is not applicable here because modular groups of the resolution CRes are reduced (not all automorphisms are present).  
A simple counterexample: consider a surface group S of genus 2 and a subgroup H that is the surface group of genus 3. Let CRes be the set of homs from H to a free group that are extendable to homs from S. Not all elements from S can be given by a formula in elements from H.  He would better used framed resolutions there himself.

(48) There is no such sentence on page 542.

(49)  We do not use ``core resolutions'' and ``penetrated core resolutions'' but 
out tight enveloping fundamental sequences have properties similar to penetrated core resolutions. When we  compare $TEnv(W_{n-1})$ and $TEnv(W_{n})$  we change these tight enveloping fundamental sequences, so they may loose some of the properties but we have reduction in ``Kurosh rank'' or in size if the Kurosh rank is the same.

(50-51) Derivatives of the misunderstanding of our notion of a tight enveloping fundamental sequence.
 
(52)  It is the phrasing that is criticized here. We say `` the sequence stabilizes'',
Sela said `` if there is no change in any of them''. We are talking about the same thing here.  

(53)   No comments.

(54) Sela does not explain why the statement is wrong.

(55)  We do not say the statement is obvious, we don't prove it referring to Theorem 11, maybe we should write down the complete proof. 

(56)  See response to comment (47) about generic families.  It worth mentioning that Sela himself does not provide a proof of  any form  of ``Generalized Merzliakov theorem''.  

In the case when QH subgroups $Q_1,\ldots ,Q_n$ of an NTQ system are finite index subgroups of larger QH subgroups $\bar Q_1,\ldots ,\bar Q_n$ of an ambient NTQ system, and the fundamental sequence is a restriction of the fundamental sequence corresponding to the larger NTQ group, we have reduced groups of automorphisms corresponding to  
$Q_1,\ldots ,Q_n$.  It is obvious that we can not always  lift solutions into $Q_1,\ldots ,Q_n$  but 
we can lift them into some subgroups $\hat Q_1,\ldots ,\hat Q_n$ such that $Q_i\leq\hat Q_i\leq \bar Q_i,$
$i=1,\ldots , n.$  Sela calls  $\hat Q_1,\ldots ,\hat Q_n$ a frame, we do not have a similar special notion.

We looked at Theorem 1.33 in \cite{Sela5} and Lemma 6 in \cite{Sela51} and found differences in their formulations as Sela suggests.  It would be nice to compare their proofs but, unfortunately, it is written that the proof of Lemma 6  is the same as the proof of Theorem 1.33 (inspite of differences in the formulations). Theorem 1.33  has a very unclear proof.

(57) Sela writes that C. Perin \cite{Perin} found  that his \cite{Sela6}, Theorem 7 describing f.g. groups elementarily equivalent to a nonabelian free group, is wrong. He claims that our \cite{KMel}, Theorem 41 is wrong too. But our theorem is stated using our concept of regular NTQ groups and is correct.  One can check which quadratic equations are regular in \cite{KMel}, Corollary 1,  and what means for a quadratic equation to have a non-commutative solution. What is means is that  regular NTQ groups are not completely identical to hyperbolic $\omega$-residually free towers, and that Sela adopted them without complete understanding.

Sela's Theorem 7 has  mistakes of different types. We will give two different examples.

a) According to \cite{Sela6}, page 8, lines -2--1, 
a fundamental group $G$ of a closed non-orientable surface of genus 3 is a hyperbolic $\omega$-residually free tower, therefore by his Theorem 7, should be elementarily equivalent to a non-abelian free group. But all solutions of the equation $x^2y^2z^2=1$ commute in a free group. Therefore, the formula $$\forall x,y,z (x^2y^2z^2=1\rightarrow [x,y]=[x,z]=[y,z]=1)$$ is true in a free group, but false in $G$.  Moreover, $G$ is not even a fully residually free group.

b) Consider two groups  $$G_1=<a,x_1.x_2,x_3|x_1^{-1}ax_1x_2^{-1}a^{-1}x_2x_3^{-1}ax_3a^{-1}=1>$$ and $$G_2=<a,x_1,x_2,x_3|x_1^{-1}ax_1x_2^{-1}ax_2x_3^{-1}ax_3a^{-3}=1>.$$

Both groups are non-elementary hyperbolic by the combination theorem (as one separated HNN extension with the base group $F(a,x_2,x_3)$ and stable letter $x_1$ along elements that are not proper powers). Both groups are obtained by amalgamating a 4-punctured sphere along its boundary components (\cite{Sela6}, page 9) to the group associated with the previous level (a free group $<a>$). There is, certainly, a retraction from $G_1$ and $G_2$ onto  $<a>$ mapping all the generators to $a$. According to Sela's definition of a hyperbolic $\omega$-residually free tower in \cite{Sela6} and Theorem 7 they are both elementarily equivalent to a free non-abelian group. 

According to our Theorem 41, $G_1$ is a  regular NTQ group and is elementarily equivalent to a non-abelian free group (indeed, there is a non-commutative solution of the equation $x_1=x_2=b, x_3=a$ in a free non-abelian group).  The group $G_2$ is not even a fully residually free group, because all solutions of the equation 
$x_1^{-1}ax_1x_2^{-1}ax_2x_3^{-1}ax_3a^{-3}=1$ in a free group commute. Therefore $G_2$ is not elementarily equivalent to a free non-abelian group.

Basically, one should not refer to Sela's Theorem 7, it is completely wrong.

\section{The state of the theory of Algebraic Geometry in free groups from our point of view}
 Tarski's  problems   have attracted the attention of algebraists and logicians for a long time. Their solution is not a single act and not the result of a single paper. In fact these problems required the creation of several new  areas of group theory and a concerted effort of a large group of mathematicians which spreads across several generations. 

Even the final push that culminated in ours and Sela's works involves many other   interesting results and ideas, some of them obtained jointly with  other people. Our proof is based on several publications of total length of 800 pages.  Fortunately, it is not a single extremely long argument which is designed just  to prove the theorem; the sequence of papers addresses seemingly separate well-known problems, some of which are as interesting as the main theorems themselves.  This work provided an impetus to  the study of equations in groups, group actions, length functions, algorithmic group theory,  limits of groups, model theory of groups and algebraic geometry of groups, to mention a few.  

The development and publication of all these results was not a polished ordered process; it involved the  creation of many new concepts, numerous adjustments, changes in exposition; large pieces of the argument had to be rewritten when a better means of explaining became available. It was a process of laying down of fundamentals of some new theories, of the rethinking of some old ones,  of making new connections between existing ones; all along the way  proving and polishing our results and changing our explanations as required.  We think that the whole subject is still under construction, and undergoing many changes, improvements,  and streamlining of the argument.  
Do we (and Sela) have to correct some inessential mistakes that occurres during the work? Of course, we do. The papers contain hundreds of pages of very specialized proofs, there are errors. 
Are there any more mistakes, gaps in the argument,  badly explained logical steps in the proofs? Probably. Sela's   Theorem 7  in \cite{Sela6} is not correct, this is also Theorem 6 in  \cite{SelaICM}. Sela's Theorem 5.12, \cite{Sela1} and the statements in Definition 5.11 \cite{Sela1} are also wrong (see, for comparison, Proposition 4.21, \cite{Gui}). But our essential arguments and proofs are sound, and most of the others have been corrected already in the subsequent papers (by us or some other mathematicians).  The sheer number  of new results based on our ideas, concepts and techniques, obtained during and after our work on Tarski's problems confirms that the argument is solid and correct. Will it be improved later? No doubt.  We proved a great number of algorithmic results, beginning with the  algorithm to find the abelian JSJ decomposition of a limit group \cite{JSJ}, and solved algorithmically many Algebraic Geometry over free groups questions (see references on pages 508-514, \cite{KMel}), finding irreducible components of finite systems of equations, analogs of elimination and parametrization theorems in classical algebraic geometry etc.

Sela's work was greatly influenced by  our results on f.g. fully residually free groups published in 1998  \cite{Irk},\cite{BMR} (he renamed them limit groups or used the terminologyy "$\omega$-residually free groups'') see below. He does not refer to our results not only on Tarski's problems but also  on algebraic geometry for a free group.

We think that looking  for small inessential errors (which have been corrected afterwards) is an exercise in futility.  It makes sense only to discuss those pieces of the argument which are not clear even now, with the intention to improve. 

An absolutely  separate issue is to address  Sela's claims that some ideas were stolen from his work.  This is just ridiculous. We neither have time nor desire to have this endless and meaningless argument about priority of the results. 

 Sela  recently complained to the AMS Committee on Professional Ethics.
The committee decided that our ``papers received appropriate refereeing, both for their mathematical content and for the issue of plagiarism, which was raised in (Sela's) complaints to the editor at the time. The referees considered these charges even as they were working through the mathematics involved. In the end, neither referee thought the plagiarism charges should stand... The referees also recommended publication of the long papers, after asking for and getting substantial revisions. In our judgment an appropriate process was carried out.''

 One  has to look at the  whole work on Tarski's problems by many people to see where the ideas really came from, why we did it this way or that way.  Our approach to Tarski's problems stems from a huge body of work that was done in Malcev's, Adjan's, and Lyndon's schools of combinatorial algebra and algorithmic group theory. Sela's approach is more topological.    Quite often  both approaches  describe the same phenomena, though from different view-points; this is  a natural development, and not a sign of theft.  On the contrary, we think this is a plus as it clarifies and enriches ideas of both. In fact, the future polished exposition should be a combination of the two.  In paper \cite{KMhyp} we made an attempt to bring both approaches closer together,  to underline the similarities, and to produce a kind of a dictionary that translates one work into the language of the other. 

This attempt possibly led to  Sela's ``report''. The objects in the two works are similar but not identical. They are not amenable to a crude direct interpretation; some of our statements would not be true if interpreted via a ``computer translation'' into his language (and vice versa)
(see, for example, response to his comment 57).  Many of his critical comments resulted from such an exact translation of our concepts into his language. Additional misrepresentations result from  not remembering that some statement was made two pages before (such as  that we only  consider fundamental sequences satisfying first and second restrictions).  He is trying to prove that at the end of our  paper \cite{KMel} we have the same constructions and procedures as he does and they are different because we ``misunderstood'' his constructions
and "messed up all that''.  He does not talk to us, and does not ask questions that would clarify what was actually meant.  

In conclusion, our results are correct and provide a viable methodology that we can build on to extend the existing results on decidability to various important classes of groups like free products of groups with decidable theory, Right Angled Artin groups and in other direction.

\section{Inaccuracies and similarities in Sela's papers}
There are  mistakes in Sela's own papers. 
Sela's original papers on Tarski had many essential mistakes (fixed in the works of other people).   By the time we found that the main
result of his preprint 3 on Diophantine Geometry was wrong, Sela already had been informed about his mistake by other mathematicians. He fixed
the paper and the published version is better (it is hard to say whether it is correct or not because many of the results don't have proofs),
but he probably still has the wrong 2001 version on his home page.

a) Sela's  \cite{Sela6}, Theorem 7 is wrong (see response to comment 57 for details).  

b) Theorem 5.12, \cite{Sela1} and the statements in Definition 5.11 \cite{Sela1} are also wrong (see, for comparison, Proposition 4.21, \cite{Gui}).

Sela's published papers have  many  statements without proofs. Even old results, which are correct,  like  `` Sela's shortening argument'',  do not  have  written proofs. Some oral tradition is maintained whereby one can only understand his papers on Tarski  after spending hours   listening  to the author's explanations. The actual text does not contain enough information to enable us to check the details.

 Sela's  papers \cite{Sela3}, \cite{Sela4}, \cite{Sela5} on Tarski are unreadable (we could not read them) and contain statements that do not make sense, for example:
 
 c)  `` defining relations ...correspond to a closed loop in the $R$-state'' \cite{Sela3}, p.26, line -13-14 (which is a rescaled Cayley graph of the free group, therefore a tree, and these elements are trivial.)
 
 d)  \cite{Sela4}, p.88, lines -6- -3 ``Clearly, if an element $y_0\in F_k$ can be extended to a shortest form specialization with respect to one of the constructed closures and this shortest form $\ldots$, then the sentence is valid for $y_0$''.  Definition 4.1 of shortest form specializations is such that not every element $y_0\in F_k$ that is a specialization of some limit group $Rlim _i(y,a)$, $i=1,\ldots , m$ (we are using Sela's notation) can be extended to a shortest form specialization with respect to one of the constructed closures. Therefore these specializations are lost and what is happening after does not make any sense.

 Example: Let $Rlim(y_1,y_2, a_1,a_2)$ be an amalgamated product
 $$<y_1,y_2>=_{[y_1,y_2]=[a_1,a_2]}<a_1,a_2>.$$ The well-separated resolution $Res(y,a)$ has two levels and copnsides with its completion and closure. The only shortest form specialization here is $y_1=a_1, y_2=a_2.$ All the other values of $y_1,y_2$ are lost.
 
e) \cite{Sela4}, page 95, lines -8- (-4). He states that we can associate the set of closures and formal solutions with the resolution CRes  and refers to Theorem 1.18 in \cite{Sela2}. (This is one of the form of the implicit function theorem) This theorem is not applicable here because modular groups of the resolution CRes are reduced (not all automorphisms are present).  

A simple counterexample: consider a surface group S of genus 2 and a subgroup H that is the surface group of genus 3. Let CRes be the set of homs from H to a free group that are extendable to homs from S. Not all elements from S can be given by a formula in elements from H. 

All  [Sela4] constructions are based on this. One also should say that Theorem 1.18 in [Sela2] is not proved, so he is confused about what can be proved and how to apply it.
 
 f) Theorem 3.8, \cite{Sela7} about a global bound on the number of `` strictly solid'' families. Not even a hint of the proof.
 
  We refer to these papers because we follow the tradition in the mathematical community to refer to other people's papers on the same subject.

Regarding Sela's claim about ``stealing'' the ideas.  First of all, we always refer to Sela's papers.  Secondy, there are, indeed, some  ``similarities''.

{\bf Some  ``similarities'':}
 
1) The idea to change Makanin-Razborov's fundamental sequences and to use the JSJ decomposition appeared in our preprint of 1999 (three years before Sela's preprints). We still needed an analog of Razborov's process to get algorithmic results. His "resolutions" and "towers" are just younger  sisters of our "fundamental sequences" and "NTQ systems."
To be more precise,
Corollary 4.4 (which is an important result)  in Sela's paper \cite{Sela1} is a
reformulation of Corollary 3 and 4
of Theorem 6 in our  paper \cite{Irk}.
The object that he calls( in section 6 of [4) omega-residually free towers
was introduced in our paper \cite{Irk} and is called a coordinate group of a
nondegenerate triangular quasi quadratic (NTQ) system. 
Claim 7.5 in his paper is Corollary 2 in our paper \cite{Irk}.

Theorems 3.2, 4.1 and 8.2 in \cite{Sela1} are corollaries of Theorem 6 \cite{Irk}.  There is no reference.  Proposition 5.1 in \cite{Sela1} follows from  \cite{Guba} and Lemma 1 in \cite{Irk} there is no reference.

2) Sela never refers to Baumslag, Myasnikov, Remeslennikov's papers (published in 1999) \cite{BMR} on algebraic geometry for a free group. The whole conceptual basis that he uses was developed there, yet he does not provide a single reference.
 
3) Sela is not mentioning at all our long paper \cite{Imp} (the preprint was available in 1999). It is not in his interest to talk about this paper.   Sela's notion of "completions, closures and formal solutions" (his paper \cite{Sela2}) similar to the same results and notions from our "Implicit function over free groups" paper.
 Theorem 1.18, the main result of \cite{Sela3} is, basically, Theorem 12 in \cite{Imp}.
Our paper was widely distributed at many conferences.  Our Implicit function theorem (=Sela's formal solutions) is the key point in the proof of Tarski's conjectures.  Sela does not even prove the results in \cite{Sela2}, he just states them, so it is hard to say whether he has any proof or not.  But the proof was already provided by us in \cite{Imp}. However he does not refer to any results of \cite{Imp}. Our approach was developed in \cite{Imp}, and  he  uses that same approach without any references.

4) The process (termination of which one has to prove) was suggested by us. It is hardly surprising that the proofs have similarities.    In addition, the structure of definable sets in a free group determines some key moments in the proof. He  always insists that the main body of our work is  paper \cite{KMel} ignoring our previous papers, because he uses the content of these previous papers without references.

\end{document}